\documentclass{article}

\usepackage{graphicx} % Required for inserting images
\usepackage{amsmath}
\usepackage{amssymb}
\usepackage[round]{natbib}
\usepackage{url}
\usepackage{algpseudocode}
\usepackage{algorithm}
\usepackage{booktabs}
\usepackage{multirow}
\usepackage{colortbl}
\usepackage{authblk}
\usepackage{xcolor}

\usepackage[top=35mm, bottom=40mm, left=35mm , right=35mm]{geometry}

\DeclareMathOperator{\argmax}{argmax}

\title{
 Parameter estimation of hidden Markov models: comparison of EM and quasi-Newton methods with a new hybrid algorithm.}

\author[1,2,3,*]{Sidonie Foulon}
\author[1,2]{Thérèse Truong}
\author[3]{Anne-Louise Leutenegger}
\author[1]{Hervé Perdry}
\affil[1]{Université Paris-Saclay, UVSQ, Inserm, CESP, 16 av PV Couturier, 94807, Villejuif, France}
\affil[2]{Institut Gustave Roussy, 114 rue Edouard-Vaillant, 94805, Villejuif, France }
\affil[3]{NeuroDiderot, Inserm, Université Paris Cité, UMR1141, 48 bd Sérurier, 75019, Paris, France}
\affil[*]{Corresponding author: Sidonie Foulon, sidonie.foulon@inserm.fr}
\date{2024}

\begin{document}

\maketitle

\section*{Abstract}

Hidden Markov Models (HMM) model a sequence of observations that are dependent on a hidden (or latent) state that follow a Markov chain. These models are widely used in diverse fields including ecology, speech recognition, and genetics.

Parameter estimation in HMM is typically performed using the Baum-Welch algorithm, a special case of the Expectation-Maximisation (EM) algorithm. While this method guarantee the convergence to a local maximum, its convergence rates is usually slow.

Alternative methods, such as the direct maximisation of the likelihood using quasi-Newton methods (such as L-BFGS-B) can offer faster convergence but can be more complicated to implement due to challenges to deal with the presence of bounds on the space of parameters.

We propose a novel hybrid algorithm, QNEM, that combines the Baum-Welch and the quasi-Newton algorithms. QNEM aims to leverage the strength of both algorithms by switching from one method to the other based on the convexity of the likelihood function.

We conducted a comparative analysis between QNEM, the Baum-Welch algorithm, an EM acceleration algorithm called SQUAREM (Varadhan, 2008, Scand J Statist), and the L-BFGS-B quasi-Newton method by applying these algorithms to four examples built on different models. We estimated the parameters of each model using the different algorithms and evaluated their performances.

Our results show that the best-performing algorithm depends on the model considered. QNEM performs well overall, always being faster or equivalent to L-BFGS-B. The Baum-Welch and SQUAREM algorithms are faster than the quasi-Newton and QNEM algorithms in certain scenarios with multiple optimum. In conclusion, QNEM offers a promising alternative to existing algorithms.

\section*{Keywords}

Hidden Markov models, Baum-Welch, quasi-Newton, L-BFGS-B, SQUAREM, optimisation, computational statistics

\section{Introduction}

Hidden Markov Models (HMMs) are widely used in time series analysis, with application across diverse fields, including ecology, speech recognition and many others. A particular interest is found in biology, especially in genetics: HMM are often used for genome annotation. For example, CpG islands identification (regions of the genome with a high frequency of cytosine-guanine dinucleotides) \citep{wu2010redefining}, can be performed through the use of HMM. It can also be used to differentiate exons (protein-coding sequences) and introns (non-coding sequences) within genes \citep{arimura2005identification}. 
. 

To compute the maximum likelihood estimators of the parameters of HMM, two approaches are naturally employed: the Baum-Welch algorithm \citep{baum1966statistical, welch2003hidden}, which is a special case of the Expectation-Maximisation (EM) algorithm \citep{dempster1977maximum} (cf section \ref{BW} below) and direct maximisation of the likelihood using a quasi-Newton method (section \ref{direct}), especially BFGS (Broyden-Fletcher-Goldfarb-Shanno) algorithm \citep{broyden1970convergence, fletcher1970new, goldfarb1970family, shanno1970conditioning}. Both methods have advantages and drawbacks. In particular, the Baum-Welch algorithm accommodates naturally the presence of bounds on the space of parameters, which can be a hurdle for the quasi-Newton, in particular when some parameters live in an open interval. On the other hand, the Baum-Welch algorithm usually converges more slowly than the quasi-Newton method.

In this paper we propose a new algorithm, QNEM, for HMM parameter estimation. This algorithm starts with one or more EM steps, and switches to a quasi-Newton algorithm, the BFGS, as soon as a criterion on local convexity of the log-likelihood is met. Quasi-Newton iterations are then performed, until convergence or until the aforementioned criterion on local convexity is no longer satisfied, in which case it switches back to the EM algorithm, etc.

We compare its performances with several algorithms: the Baum-Welch algorithm, an accelerated EM algorithm called SQUAREM (section \ref{SQUAREM}) and the L-BFGS-B (Limited-memory BFGS with Bound constraints) algorithm, the ``off-the-shelf'' quasi-Newton method optimized for problems with bounded parameters.

We first review the theory of HMMs, and then expose the four algorithms of interest as well as the examples on which the comparisons are performed. Finally, we present the results obtained and discuss them.

    \section{Hidden Markov Models}

\begin{figure}[!htbp]
    \centering
    \includegraphics[width=1\linewidth]{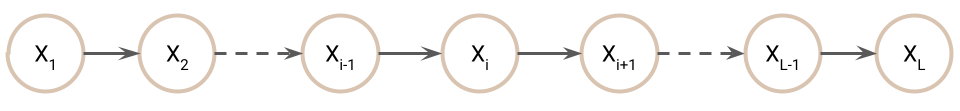}
    \caption{Representation of a Markov chain. The $X_i$ are the observations.}
    \label{fig:MaCh}
\end{figure}

A sequence of random variables measured at successive moments $(X_i)_{i \in \mathbb{N}}$ is a Markov chain if it satisfies the Markov property: for any time $i$, the distribution of $X_{i+1}$ conditional to $(X_j)_{j\leq i}$ is equal to the distribution of $X_{i+1}$ conditional to the single value $X_i$. 
The model is classically illustrated by a directed acyclic graph (DAG) in Figure \ref{fig:MaCh}.

\begin{figure}[!htbp]
    \centering
    \includegraphics[width=1\linewidth]{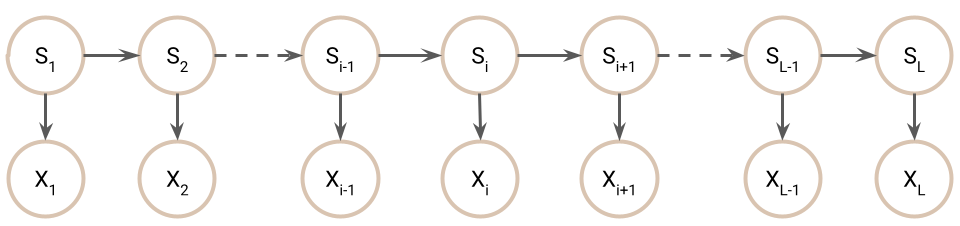}
    \caption{Representation of an Hidden Markov chain. The $X_i$ represents the observations and the $S_i$, the hidden states.}
    \label{fig:HMM}
\end{figure}

One of the property of a Markov chain is the homogeneity of the sequence. A sequence is homogeneous if the transition probabilities $\mathbb{P}(X_{i+1} | X_i)$ are identical in every point of the sequence. This property is sometimes contradicted by a given sequence of observations, making the Markov model inappropriate. A more complex model, in which there is an additional hidden layer $(S_i)_{i \in \mathbb{N}}$ which is a Markov chain, and the distribution of $X_i$ is assumed to depend only on the ``hidden state'' $S_i$, can be more appropriate. This model is called Hidden Markov Model (HMM) and is illustrated by the DAG in Figure \ref{fig:HMM}.

In some situations, the Markov model is known to be adequate, but it can be impossible to directly access the Markov chain, but only variables depending on the current state of the Markov chain. In this case, the HMM provides a natural framework for analysis. HMMs allow to retrieve the states $S_i$ of the latent Markov chain, using the observations $X_i$ whose distribution depends on the hidden states.

        \subsection{Direct maximisation of the likelihood with a quasi-Newton type method} \label{direct}

Direct maximisation of the likelihood (DML) with quasi-Newton methods is a possible way to estimate the HMM parameters. 

Quasi-Newton methods are particularly suited for maximisation of differentiable scalar functions, such as likelihood functions. They have the advantage of avoiding the calculation of the Hessian matrix by approximating the Hessian matrix (or its inverse) at each iteration. 
The off-the-shelf methods are BFGS or L-BFGS-B when some parameters are constrained, standing respectively for Broyden-Fletcher-Goldfarb-Shanno \citep{broyden1970convergence, fletcher1970new, goldfarb1970family, shanno1970conditioning}, and for Limited-memory BFGS with Bound constraints. 

This method relies on a forward step, computing probabilities of the hidden states, either conditionally to the $X_i$, or jointly with the $X_i$. These probabilities allow to compute the likelihood. The most common choice for the DML is using joint probabilities going forward \citep{thompson2000statistical, welch2003hidden, leutenegger2003estimation}. A log-transformation is often employed to get rid of the numerical underflow implied by the joint probabilities \citep{zucchini2009hidden}. In this paper, we chose to use the conditional version of the algorithm for the forward step (cf below section \ref{Estep}), which avoids the numerical underflow and therefore going through the log-transformation.

        \subsection{Baum-Welch algorithm, an EM algorithm} \label{BW}

The Baum-Welch algorithm \citep{baum1966statistical, welch2003hidden} was developed to estimate the HMM parameters.

The respective contributions of Baum and Welch in the development of the Baum-Welch algorithm is explained in Welch's Shannon lecture, in 2003 \citep{welch2003hidden}. After working separately in the calculation of a posteriori probabilities of the hidden states, Baum and Welch worked together on the re-estimation of the parameters from these a posteriori probabilities. 

The Baum-Welch algorithm was later recognised as a particular case of the more general Expectation-Maximisation (EM) algorithm formalised by Dempster, Laird and Rubin in 1977 \citep{dempster1977maximum}. 

The Baum-Welch algorithm relies on a forward-backward algorithm for calculating the a posteriori probabilities of the hidden states (E-step). This algorithm is constructed with a forward step as sketched before. We chose to use the forward algorithm based on conditional probabilities, which is more commonly used for the Baum-Welch algorithm. It is also possible to use the forward algorithm based on the joint probabilities, together with an adequate version of the backward step (cf \citep{thompson2000statistical}, chapter 6).
Then the parameters are re-estimated from the a posteriori probabilities (M-step). This step depends on the model considered.

\section{Material and methods}

\newcommand{\T}{T_\theta}
\newcommand{\E}{E_\theta}

For the following, let $ \textbf{X} = (X_1, ..., X_L) = X_1^L$ the observations and $\textbf{S} = (S_1, ..., S_L) = S_1^L$ the hidden states along the chain of length $L$.

Let $\T(s, t) = \mathbb{P}_\theta(S_{i+1} = t| S_i = s)$ be the transition probability from state $s$ to state $t$ and $\E(x,s) = \mathbb{P}_\theta(X_i = x| S_i = s)$ the emission probability of observation $x$ from state $s$.

We consider four maximum likelihood methods to estimate the parameter $\theta$.

\subsection{Direct maximisation of the likelihood}\label{direct2}

As seen in for example Zucchini \citep{zucchini2009hidden} and Turner \citep{turner2008direct}, a variant of the forward algorithm based on joint probabilities can be used to compute the likelihood. If $a_i(s) = \mathbb{P}(S_i = s, X_1^{i-1} = x_1^{i-1})$ and $b_i(s) = \mathbb{P}_\theta(S_i = s,X_1^{i} = x_1^{i})$, the quantities $a_i(s)$ and $b_i(s)$ can be computed recursively as follows:

\begin{enumerate}
    \item Initialisation : \\
    \begin{equation}\label{eq:fw a}
        a_1(s) = \pi(s)
    \end{equation}
    \begin{equation}
        b_1(s) = a_1(s) \E(x_1, s)
    \end{equation}
    with $\pi(s)$ the initial distribution of state $s$.
    \item Recursion : \\
    \begin{equation}
        a_i(s) = \sum_t b_{i-1}(t) \T(t, s)
    \end{equation}
    \begin{equation}
        b_i(s) = a_i(s) \E(x_i, s)
    \end{equation}
    for $i \in {2, ..., L}$
\end{enumerate}

Finally the likelihood can be computed as
    \begin{equation} \label{eq:lik}
        L(\theta ; \textbf{X}) = \mathbb{P}_\theta(X_1^L = x_1^L) = \sum_s a_L(s)
    \end{equation}

As the values of the $a_i(s)$ and $b_i(s)$ become very small when $i$ increases, this algorithm is prone to numerical underflow. The usual solution is to compute the values on the logarithmic scale \citep{zucchini2009hidden}, which involves rewriting the above equations accordingly. In this paper, we prefer to use the version of the forward algorithm based on conditional probabilities.

This forward algorithm computes classically two conditional probabilities alternatively: 
\begin{equation} \label{eq: fw alpha}
    \alpha_i(s) = \mathbb{P}_{\theta}(S_i=s | X_1^{i-1} = x_1^{i-1})
\end{equation}
the ``forecast'' probability and
\begin{equation} \label{eq: fw beta}
    \beta_i(s) = \mathbb{P}_{\theta}(S_i=s | X_1^{i} = x_1^{i})
\end{equation}
the ``filtering'' probability. 

To compute the log-likelihood we need to compute also for each $i$
$\gamma_i = \log\left(\mathbb{P}_\theta(X_1^i = x_1^i) \right)$. These quantities can be computed recursively as follows:

\begin{enumerate}
    \item Initialisation : \\

    \begin{equation}\label{eq:init alpha}
        \alpha_1(s) = \pi_s 
    \end{equation} 
    \begin{equation}
        \beta_1(s) = \frac{\E(x_1,s) \times \alpha_1(s)}{\sum_i \E(x_1,i) \times \alpha_1(i)}
    \end{equation} 
    \begin{equation}
         \gamma_1 = \log\left(\sum_s \alpha_1(s) \E(x_1, s)\right).
    \end{equation}

    \item Recursion : \\
 
    \begin{equation}
        \alpha_i(s) = \sum_t \T(t,s) \times \beta_{i-1}(t)
    \end{equation} 
    \begin{equation}
        \beta_i(s) = \frac{\E(x_i,s) \times \alpha_i(s)}{\sum_t \E(x_i,t) \times \alpha_i(t)}
    \end{equation} 
    \begin{equation} \label{eq:gamma}
        \gamma_i = \gamma_{i-1} + \log\left(\sum_s \E(x_i, s) \alpha_i(s)\right)
    \end{equation}
    for $i \in {2, ..., L}$
\end{enumerate}
Finally, the log-likelihood $\ell(\theta ; \mathbf{X})$ is equal to $\gamma_L$. The likelihood is still computed on the logarithmic scale, but the $\alpha_i(s)$ and $\beta_i(s)$ can be computed without rewriting all the equations on the logarithmic scale, which make the computation simpler.

The L-BFGS-B algorithm is then performed on the likelihood computed in the last step using the R \texttt{optim} function. The gradient of the likelihood along $\theta$ can be computed using automatic differentiation (our implementation relies on the R package salad \citep{perdry2024salad}). 

    \subsection{Baum-Welch algorithm} 

The Baum-Welch EM algorithm is composed by two steps:
\begin{itemize}
    \item E step: knowing $\theta$, compute the a posteriori probabilities of the hidden states along the chain (forward-backward algorithm)
    \item M step: re-estimate $\theta$ using these probabilities
\end{itemize}
The Baum-Welch algorithm alternates these 2 steps until convergence. As mentioned in section \ref{BW}, we will use here the forward-backward algorithm presented in \citep{zucchini2009hidden} or \citep{delmas2006modeles}. 

        \subsubsection{E step : Forward-Backward algorithm} \label{Estep}

The forward-backward algorithm runs through the observations twice: the first time forward, ie. from first to last observations ; the second backward, ie. from last to first observation. 

After computing the forward quantities shown in equations (\ref{eq: fw alpha}) and (\ref{eq: fw beta}), the following ``backward quantities'' will be computed:
\begin{equation}
    \delta_i(s,t) = \mathbb{P}_{\theta}(S_{i-1}=s, S_i=t | X_1^{L} = x_1^{L})
\end{equation}
the ``smoothing" probability and
\begin{equation}
    \varphi_i(s) = \mathbb{P}_{\theta}(S_i=s | X_1^{L} = x_1^{L})
\end{equation}
the ``marginal" probability. 

\begin{enumerate}
    \item Initialisation: 

    \begin{equation}
        \varphi_L(s) = \beta_L(s) 
    \end{equation}
    \item Recursion:

    \begin{equation}
        \delta_i(s,t) = \T(s,t) \times \frac{\beta_{i-1}(s)}{\alpha_i(t)} \times \varphi_i(t)
    \end{equation}
    \begin{equation}
        \varphi_{i-1}(s) = \sum_t \delta_i(s,t)
    \end{equation}
for $i \in {L, ..., 2}$.
\end{enumerate}

At the end of the forward-backward, we have computed the probabilities $\varphi_i(s)$ of the hidden states, and the probability $\delta_i(s,t)$ of two consecutive hidden states, given all the observations. These probabilities will be used in the M step.\\

It is worth to note that the backward algorithm has a higher computational cost than the forward algorithm, due to the necessity to calculate the $\delta_i(s,t)$ for all pairs of states $s$ and $t$, at each time point $i$.
As mentioned previously, there exists a variant of the backward algorithm based the forward quantities $a_i$ and $b_i$ involving joint probabilities, as seen in \citep{thompson2000statistical}.

        \subsubsection{M step}

The aim in M step is to re-estimate the parameter $\theta$ that will be used in the following E step.

Note $\theta^{(k)}$ the parameter at the $k$-th iteration. For each iteration $w$ of the Baum-Welch algorithm, in M step we are going to estimate the parameter $\theta^{(k+1)}$ for the next iteration.

The complete likelihood of $\theta$, ie the likelihood assuming the hidden states are observed, is: 
\begin{equation}
    \begin{split}
        L(\theta ; \textbf{S} = \textbf{s}, \textbf{X} = \textbf{x}) =  \mathbb{P}_{\theta}(S_1 = s_1) & \cdot \prod_{i=2}^L \mathbb{P}_{\theta}(S_i = s_i | S_{i-1} = s_{i-1}) \\
        & \cdot \prod_{i=1}^L \mathbb{P}_{\theta}(X_i = x_i | S_i = s_i)
    \end{split}
\end{equation}

Now the log-likelihood of $\theta$ is: 
\begin{equation}
    \begin{split}
        \ell(\theta ; \textbf{S} = \textbf{s}, \textbf{X} = \textbf{x}) =  \log \mathbb{P}_{\theta}(S_1 = s_1) & + \sum_{i=2}^L \log \mathbb{P}_{\theta}(S_i = s_i | S_{i-1} = s_{i-1}) \\
        & + \sum_{i=1}^L \log \mathbb{P}_{\theta}(X_i = x_i | S_i = s_i)
    \end{split}
\end{equation}

The M step relies on the expected value of the log-likelihood: 
\begin{equation}
    \begin{split}
        Q(\theta ; \theta^{(k)}) & = \mathbb{E}( \ell(\theta ; \textbf{S} = \textbf{s}, \textbf{X} = \textbf{x}) |  \theta^{(k)}) \\
        &= \sum_s \left( \log \mathbb{P}_{\theta}(S_1 = s) \cdot \varphi_1(s) \right) + \sum_{i=2}^L\left[ \sum_{s,t} \left( \log \mathbb{P}_{\theta}(S_i = t | S_{i-1} = s) \cdot \delta_i(s,t) \right) \right] \\
        & + \sum_{i=1}^L \left[ \sum_s \left( \log \mathbb{P}_{\theta}(X_i = x_i | S_i = s) \cdot \varphi_i(s) \right) \right]
    \end{split}
\end{equation}

The M step ends with finding the parameters that maximise $Q(\theta ; \theta^{(k)})$ : 
\begin{equation}
    \theta^{(k+1)} = \argmax_{\theta} Q(\theta ; \theta^{(k)}).
\end{equation}

In concrete examples, the $\delta_i(s,t)$ will be used to compute the transition parameters of the Markov chain, and the $\phi_i(s)$ the parameters relatives to the emission probabilities.

    \subsection{EM acceleration} \label{SQUAREM}

The SQUAREM method was presented in Varadhan, 2008 \citep{varadhan2008simple}. It is a method created to accelerate the EM algorithm. This method seems interesting when we aim to apply EM to high-dimensional data or data based on complex statistical models. In a nutshell, the SQUAREM makes two EM iterations, and uses these values to extrapolate the trajectory of the algorithm, and estimate a future point several iterations ahead. It then iterates again from this extrapolated point. The constraints on the parameters, if any, are simply dealt with by constraining the extrapolated point to be in the feasible space. We applied this method to our examples to test its performances as well.

    \subsection{QNEM: a mix of the two algorithms}
    
We propose to mix the Baum-Welch algorithm and the quasi-Newton algorithm BFGS as presented in the algorithm \ref{alg:qnem}.  

\begin{algorithm}[htb]
    \caption{QNEM algorithm}\label{alg:qnem}
    \begin{enumerate}\addtocounter{enumi}{-1}
        \item   Input : 
                $\theta^{(0)}$ = starting point. \\
                Set $k = 0$ and $H_0 = \mathbb{I}$ .
               
        \item   Compute 
        \[
            \theta^{(k+1)} = BW(\theta^{(k)})
        \]

        \item   Let $s_k = \theta^{(k+1)} - \theta^{(k)}$ and $y_k = \nabla f_{k+1} - \nabla f_k$. \\
                If the curvature condition $s'_k \cdot y_k > 0$ is not satisfied, go to 1.

        \item   Update $H_k$ using the BFGS formula
        \[
            H_{k+1} = H_k - \frac{H_k \cdot y_k \cdot y'_k \cdot H_k}{y'_k \cdot H_k \cdot y_k} + \frac{s_k \cdot s'_k}{y'_k \cdot s_k}        
        \]\\
        Let $k = k+1$.

        \item   Let $p_k = -H_k \cdot \nabla f_k$ the search direction. \\
                Perform a backtracking search for $\alpha$ such that Armijo condition is satisfied
        \[
            f(\theta_k + \alpha \cdot p_k) \leq f(\theta_k + c \cdot \alpha (\nabla f_k)' p_k)
        \]

        \item   Let $\theta^{(k+1)} = \theta^{(k)} + \alpha \cdot p_k$.\\
                If the curvature condition is not satisfied, reinitialise $H_k = \mathbb{I}$ and go back to 1.\\
                Else go to 3.
    \end{enumerate}
\end{algorithm}

If the current value of the parameter is $\theta^{(k)}$, the BFGS algorithm finds $\theta^{(k+1)}$ by performing a linear search in a search direction $p_k = -H_k \nabla f_k$ where $H_k$ is the current approximation $H_k$ of the inverse of the Hessian of the objective function, and $\nabla f_k$ is its gradient in $\theta^{(k)}$. The $H_k$ is updated using the formula
        \[
            H_{k+1} = H_k - \frac{H_k \cdot y_k \cdot y'_k \cdot H_k}{y'_k \cdot H_k \cdot y_k} + \frac{s_k \cdot s'_k}{y'_k \cdot s_k}        
        \]
where $s_k = \theta^{(k+1)} - \theta^{(k)}$ and $y_k = \nabla f_{k+1} - \nabla f_k$. If $H_k$ is positive definite, $H_{k+1}$ will be positive definite as soon as the scalar product of $s_k$ and $y_k$ is positive: 
\begin{equation} \label{eq:curv cond}
    s'_k \cdot y_k > 0
\end{equation}    
This is the curvature condition. If the function is strongly convex it always holds after a linear search which successfully decreases the value of the objective function. In the general case, it may not hold, leading to the use of complex linear search methods to avoid this situation, such as the Nocedal line search, see \citep{nocedal2006numerical} chapter 3. 

The QNEM starts by performing Baum-Welch algorithm steps, until the curvature condition (equation \ref{eq:curv cond}) is met; it then switches to BFGS algorithm. Our BFGS iterations use a simple backtracking search, until the Armijo condition is met. This condition ensures that the objective function decreased.

As long as the curvature condition holds, BFGS iterations are continued. If at some point it does not hold, the QNEM falls back to the Baum-Welch algorithm, until the curvature condition is met and the quasi-Newton iterations are performed again, until convergence.
The box constraint are simply dealt with by projecting the gradient on the box when $\theta$ is on the border.

    \subsection{Stopping criterion}

For the four algorithms, we choose to use the same stopping criterion to decide when the convergence point is reached. The criterion is the classic relative convergence tolerance (reltol), which is based on the likelihood of the current iteration and the previous one. At the $k$-th iteration, the algorithm stops if

\begin{center}
    $ \frac{\lvert \ell^{(k-1)} - \ell^{(k)} \rvert}{\lvert \ell^{(k-1)} \rvert + reltol} < reltol$
\end{center}

where $\ell^{(k)}$ is the likelihood of the $k$-th iteration and $reltol$ is a small constant, depending on the machine precision; in our experiment it is $1.49 \cdot 10^{-8}$. 
    
    \subsection{Examples}
        \subsubsection{Umbrella example}\label{methods um}

We will first use a simple example based on a fictitious situation. Every day of the year, in a secret underground installation, a security guard observes whether the director comes in with an umbrella or not \citep{gao_inference_2021}. During 56 days, the security guard keep track of his observations on his note pad as follows: U if the director carries an umbrella; N otherwise. 
The guards aims to predict the daily weather conditions (hidden states) based on the umbrella observations: rainy (noted R) or dry (noted D) day. Here the hidden state is the weather. Since the security guard is in a underground installation, he can not observe directly the weather: he has to rely only on his observations of the umbrella status to determine the weather. The parameters to estimate are \textit{a} the probability of weather state transition and \textit{b} the probability of error in the umbrella status (not carrying an umbrella on rainy day or carrying an umbrella on dry day). The transition and emission matrix are presented respectively in Tables \ref{tab:trans_u} and \ref{tab:emiss_u}.

\begin{table}[!htbp]
    \centering
    \begin{tabular}{ccc}\toprule
        Transition matrix & $S_i$ = D & $S_i$ = R\\ \midrule
        $S_{i-1}$ = D & 1-a & a\\
        $S_{i-1}$ = R & a & 1-a\\ \bottomrule
    \end{tabular}
    \caption{Transition matrix for the umbrella example}
    \label{tab:trans_u}

    \bigskip
    \begin{tabular}{ccc} \toprule
        Emission matrix & $S_i$ = D & $S_i$ = R\\ \midrule
        $X_i$ = N & 1-b & b\\
        $X_i$ = U & b & 1-b\\ \bottomrule
    \end{tabular}
    \caption{Emission matrix for the umbrella example}
    \label{tab:emiss_u}
\end{table}

To apply the quasi-Newton algorithm (L-BGFS-B) to this example, it was necessary to bound the values of $a$ and $b$ to the intervals $[0.01, 0.99]$. Using the interval $[0, 1]$ would be preferable, but the likelihood is not always defined on the boundary, which leads to a failure of the algorithm. The algorithm can not accommodate open intervals such as $(0,1)$.

        \subsubsection{Old Faithful geyser example} \label{methods gey}

Following an example presented in Zucchini \citep{zucchini2009hidden}, we work with the Old Faithful geyser data \citep{azzalini1990look}, available in the R package \texttt{MASS}, considering the 272 duration of the eruptions. 
We make two uses of the data. First, we dichotomise at three minutes of eruption to divide the observations in two groups "Duration $<$ 3min" and "Duration $\ge$ 3 min" (noted as "Dinf3" and "Dsup3" respectively) to work in a discrete observations case. This is similar to what is done in \citep{zucchini2009hidden}. 
Second, we work with the continuous values, assuming they are drawn in Gaussian distributions whose parameters depend on the hidden states.
In the sequence of dichotomised values, it appears that each short eruption was followed by a long eruption, but there are sequences of long eruptions. We decided to assume that the hidden Markov chain has three states: "short", "long" and "steady long", noted respectively "S", "L" and "Sl". 
The parameters to estimate in the dichotomised model are \textit{a} the probability to change from state "short" to "steady long", \textit{b} the probability to stay in "steady long" state, \textit{c}, \textit{d} and \textit{e} the probability to observe a "Duration $\ge$ 3 min" in respectively "short", "long" and "steady long" states. This emission matrix is presented in Table \ref{tab:emiss_gd}.
For the continuous model, the parameters to estimate are the same \textit{a} and \textit{b} and six other parameters $\mu_s$, $\mu_l$ and $\mu_{sl}$, the mean of the Gaussian distributions in which duration of eruptions are drawn, and $\sigma_s$, $\sigma_l$ and $\sigma_{sl}$, their standard deviation. 
The transition matrix, common to both the dichotomised and continuous models is presented in Table \ref{tab:trans_g}

\begin{table}[!htbp]
    \centering
    \begin{tabular}{cccc} \toprule
        Transition matrix & $S_i$ = S & $S_i$ = L & $S_i$ = Sl\\ \midrule
        $S_{i-1}$ = S & 0 & 1-a & a\\
        $S_{i-1}$ = L & 1 & 0 & 0\\
        $S_{i-1}$ = Sl & 1-b & 0 & b\\ \bottomrule
    \end{tabular}
    \caption{Transition matrix for the geyser example}
    \label{tab:trans_g}

    \bigskip

    \centering
    \begin{tabular}{cccc} \toprule
        Emission matrix & $S_i$ = S & $S_i$ = L & $S_i$ = Sl\\ \midrule
        $X_i$ = Dinf3 & 1-c & 1-d & 1-e\\
        $X_i$ = Dsup3 & c & d & e\\ \bottomrule
    \end{tabular}
    \caption{Emission matrix for the dichotomised version of the geyser example}
    \label{tab:emiss_gd}
\end{table}

Again, for the quasi-Newton L-BFGS-B algorithm, we had to constraint the probabilities $a$, $b$, $c$, $d$, and $e$ to the interval $[0.01, 0.99]$, and to set a lower bound of $0.01$ for  $\sigma_s$, $\sigma_l$, and $\sigma_{sl}$.
        
        \subsubsection{HBD segments example}\label{methods hbd}
        
\begin{figure}[!htbp]
    \centering
    \includegraphics[width=0.6\linewidth]{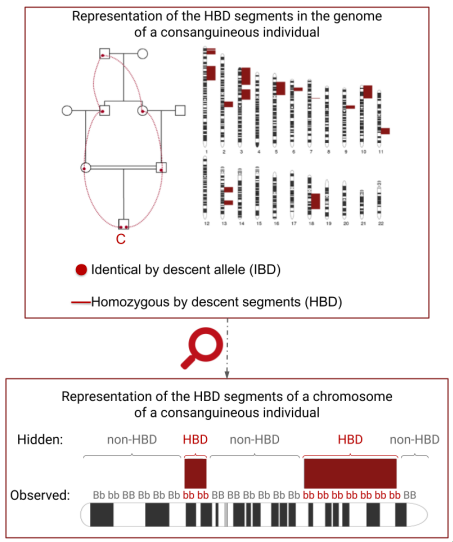}
    \caption{This figure shows on the upper left a representation of the genealogy of a consanguineous individual C and on the upper right his/her genome. In red are represented the Homozygous By Descent (HBD) segments. On the lower part is represented a chromosome with the information about the hidden states (2 possible values : HBD or non-HBD) and the observations (genotypes with B the reference allele and b the alternative allele).}
    \label{fig:cons}
\end{figure}

Consanguineous individuals are offspring of related parents. In the absence of pedigree information, consanguineous individuals can be identified through their genome which carries Homozygous-By-Descent (HBD) segments. HBD segments are regions of the genome were the two copies of the genome are identical and come from a common ancestor of the parents. These segments are characteristic of consanguineous individuals. Reconstructing these regions allow to determine whether an individual is consanguineous and his/her degree of consanguinity. The coefficient of consanguinity \textit{f} can be defined as the expected proportion of HBD genome. 
HMM are used to compute for each variant of the genome, its probability to be in a HBD segment. These probabilities can be used to reconstruct the HBD segments, but this requires to first estimate the parameters of the HMM. The parameters to estimate are the coefficient of consanguinity \textit{f} and another parameter \textit{a} related to the mean length of HBD segments \citep{leutenegger2003estimation, gazal2014fsuite}.
cf. Figure \ref{fig:cons}

To test the algorithms, we will generate the genome of an individual from a given set of parameters $\theta = (f = 0.0625, a = 0.064)$, which corresponds to the expected values for an offsprings of first cousins. The data were created from an inhouse R package named \texttt{Mozza} (available at \url{https://github.com/genostats/Mozza}, still in development). We then rounded the distances between positions to 0.1 centiMorgan (cM). For each of these new positions, we randomly selected only one observation, to create submaps of positions. This corresponds to a total of 1050 observations.

The transition and emission probabilities are presented in Tables \ref{tab:trans_h} and \ref{tab:emiss_h}. We note $d$ the distance between positions (fixed here at 0.1 cM), $pA$ and $pa$ the allele frequency of the reference and alternative allele respectively and $\epsilon$ the genotyping or mutation error. 

\begin{table}[!htbp]
    \centering
    \begin{tabular}{ccc} \toprule
        Transition matrix & $S_i$ = nHBD & $S_i$ = HBD\\ \midrule
        $S_{i-1}$ = nHBD & $(1-\exp^{-ad})(1-f) + \exp^{-ad}$ & $(1-\exp^{-ad})f$\\
        $S_{i-1}$ = HBD & $(1-\exp^{-ad})(1-f)$ & $(1-\exp^{-ad})f + \exp^{-ad}$\\ \bottomrule
    \end{tabular}
    \caption{Transition matrix for the HBD example}
    \label{tab:trans_h}

    \bigskip

    \centering
    \begin{tabular}{ccc} \toprule
        Emission matrix & $S_i$ = nHBD & $S_i$ = HBD\\ \midrule
        $X_i$ = AA & $pA^2$ & $(1-\epsilon)pA + \epsilon \cdot pA^2$\\
        $X_i$ = Aa & $2\cdot pA \cdot pa$ & $2\epsilon pA \cdot pa$\\
        $X_i$ = aa & $pa^2$ & $(1-\epsilon)pa + \epsilon \cdot pa^2$\\
        $X_i$ = NA & 1 & 1\\ \bottomrule
    \end{tabular}
    \caption{Emission matrix for the HBD example}
    \label{tab:emiss_h}
\end{table}

In this example, the probability $f$ was constrained to the interval $[0.01, 0.99]$ and the lower-bound was $0.01$ for the parameter $a$ when applying the quasi-Newton algorithm.

    \subsection{Comparison of the algorithms}
We compare the likelihood maximisation methods by comparing 
\begin{itemize}
    \item the number of iteration to convergence
    \item the number of forwards step runs
    \item the number of backwards step runs
    \item the running time
    \item the percentage of runs that failed to converge
    \item the running time of forward step
    \item the running time of backward step
    
\end{itemize}

We also studied the final estimations of the parameters and calculated the percentage of appearance of each corresponding likelihoods.

For all four examples, the tests were made on 1000 uniform draws within a wide possibility of values for the initial $\theta$, some of them being unlikely.

\section{Results}

    \subsection{Umbrella example}\label{results um}

The results are presented in Table \ref{tab:res u}. 

\begin{table}[!htbp]
    \hbox to \hsize{\hss
    \begin{tabular}{ccccccccccc} \toprule
         Umbrella example & \multicolumn{6}{c}{Nb of iterations} & \multicolumn{2}{c}{Mean nb of steps} & Mean & Percent of\\ \cmidrule(lr){2-7} \cmidrule(lr){8-9}
          & Min & Q1 & Q2 & Mean & Q3 & Max & Fw & Bw & time (s) & convergence\\ \midrule

        qN & 4 & 13 & 16 & 15.77 & 18 & 37 & 15.77 & 0 & 0.03 & 100 \\
        BW & 3 & 18 & 24 & 27.95 & 31 & 338 & 27.95 & 27.95 &  0.06 & 100 \\
        SQUAREM & 3 & 20 & 24 & 26.76 & 29 & 179 & 17.70 & 17.67 & 0.04 & 100 \\
        QNEM & 2 & 8 & 11 & 11.32 & 13 & 118 & 12.63* & 2.80 & 0.03 & 100 \\ \bottomrule
        
    \end{tabular}\hss}
    \centering
    \caption{This table recapitulates for each algorithm the number of iterations, the mean number of steps, the mean running time and the percentage of convergence in the umbrella example. The mean running times are in seconds. For all algorithms, each backward step takes $\sim$ 1.7ms, but the forward step for the Baum-Welch (BW) and the SQUAREM algorithms is $\sim$ 0.5ms only while it is $\sim$ 1.1ms for the quasi-Newton (qN) and the QNEM algorithms, due to the computation of the derivatives. * = The number of forward steps is higher than the number of iterations because the backtracking is not counted in the number of iterations.}
    \label{tab:res u}
\end{table}

We can see that, overall, the QNEM algorithm achieves convergence with the lowest mean number of iterations, averaging 11.32 across all runs. However, the quasi-Newton algorithm has the lowest maximum number of iterations (37 iterations vs. 118-338 iterations for the other three algorithms). The number of iterations of the Baum-Welch and the SQUAREM algorithms are comparable and systematically higher than those of the QNEM and the quasi-Newton algorithms.

Apart from the total number of iterations, it is interesting to look at the mean number of forward and backward steps used by the four algorithms. The quasi-Newton algorithm does not perform any backward steps. For the quasi-Newton and Baum-Welch algorithms, the number of forward steps corresponds to the mean number of iterations on the 1000 draws of the initial starting point. The Baum-Welch algorithm compute the same number of forward and backward steps. For the SQUAREM algorithm, the mean number of forward and backward steps is slightly different, with 17.70 forward steps and 17.67 backward steps. Finally, the QNEM algorithm executes an average of 12.63 forward steps and 2.80 backward steps, where the backward steps correspond to the calls to the Baum-Welch algorithm. We note that a forward step takes about 1.5ms for the Baum-Welch and the SQUAREM algorithms and around 11.5ms for the other two, which is longer because of the additional computation of the derivatives. A backward step takes about 8ms for all algorithms. 
The total mean time spent in forward and backward steps for a run is almost the same for the quasi-Newton algorithm and the QNEM ($\sim$ 18ms). It is higher for the SQUAREM, and even longer for the Baum-Welch algorithm ($\sim$ 39ms and 61.5ms respectively). This trend is in accordance with the total mean execution time, where the Baum-Welch algorithm is the slowest one and the quasi-Newton and the QNEM algorithms are the fastest.
\\

\begin{table}[!htbp]
    \hbox to \hsize{\hss
    \begin{tabular}{ccccccc} \toprule
         Likelihood & \multicolumn{2}{c}{Parameters} & \multicolumn{4}{c}{Percent of runs}\\ \cmidrule(lr){2-3} \cmidrule(lr){4-7}
          & a & b & qN & BW & SQUAREM & QNEM\\ \midrule
         \multirow{2}{*}{33.3} & \multirow{2}{*}{0.06} & 0.18 & \multirow{2}{*}{89.9} & \multirow{2}{*}{70.1} & \multirow{2}{*}{69.9} & \multirow{2}{*}{70.1} \\
          & & 0.82 & & & & \\ \arrayrulecolor{gray} \midrule
         33.8 & ND & 0.5 & 10.1 & 29.9 & 30.1 & 29.3 \\ \midrule
         Other & & & 0 & 0 & 0 & 0.6 \\ \arrayrulecolor{black} \bottomrule
        
    \end{tabular}\hss}
    \centering
    \caption{This table shows the different convergence points reached with the umbrella example and their percent. ND = Non determined.}
    \label{tab:conv u}
\end{table}

\begin{figure}[!htbp]
    \centering
    \includegraphics[width=0.6\textwidth]{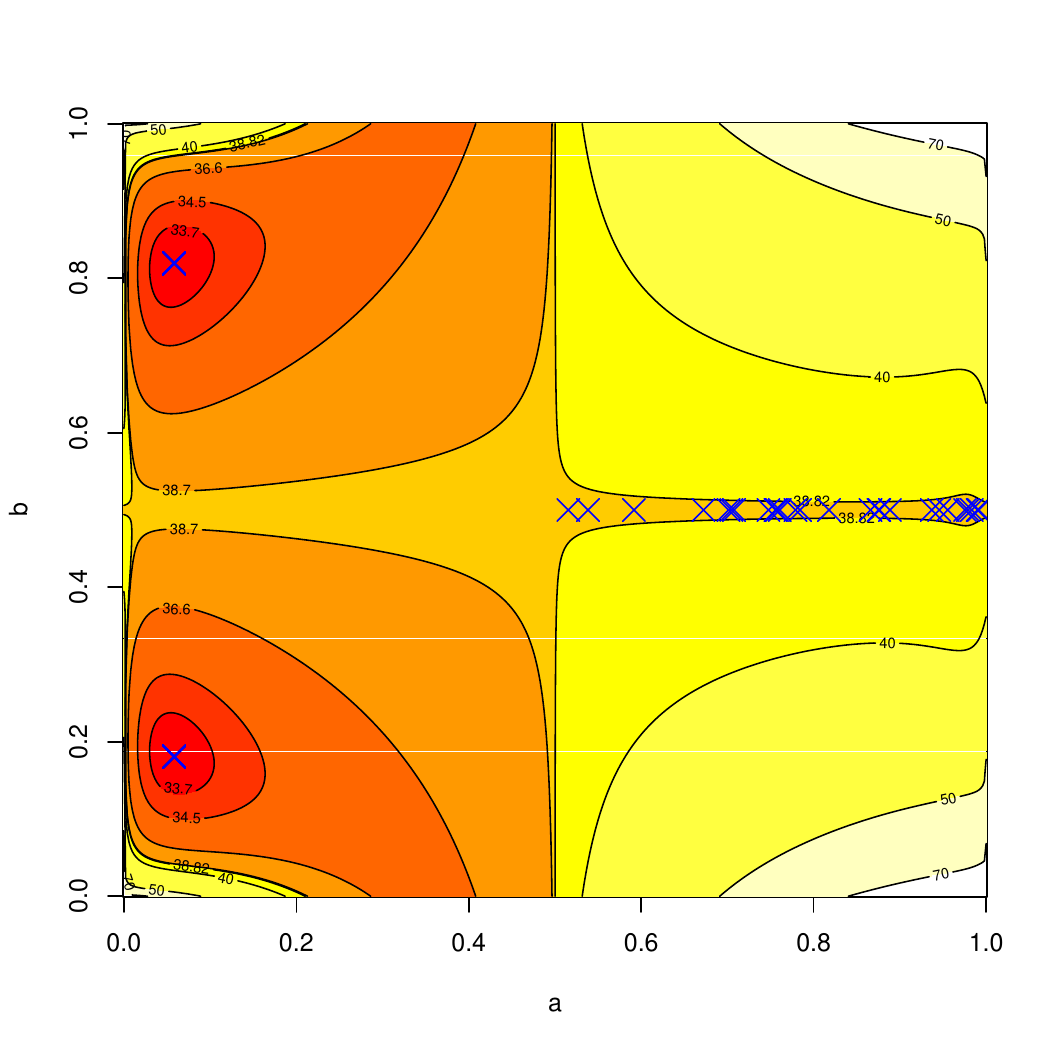}
    \caption{The likelihood with 100 convergence points.}
    \label{fig:likelihood}
\end{figure}

We also examined the final values obtained for $\theta$. They can be divided into three subsets, corresponding to three different regions of the likelihood as shown in Figure \ref{fig:likelihood}. Two of these regions are points with $a = 0.06$ and $b = 0.18$ or $b = 0.82$, corresponding to the highest likelihood value (negative log-likelihood equal to $33.3$). The third region is the segment $b = 0.5$ and $a \in [0.5, 1]$ which corresponds to a crest line (negative log-likelihood equal to $33.8$). The parameters and the corresponding proportion of runs converging to these parameters for each algorithm are presented in Table \ref{tab:conv u}. The $\theta$ with the best likelihood is reached almost 90$\%$ of the time with the quasi-Newton algorithm but only 70$\%$ of the time for the other three algorithms. 

    \subsection{Old Faithful geyser example}

We will first present the results obtained with the dichotomised data, then with the continuous values.

        \subsubsection{Dichotomised values}\label{results gd}
        
As previously done, we tested the algorithms on 1000 random initial sets of parameters. The results are presented in Table \ref{tab:res gd}.

\begin{table}[!htbp]
    \hbox to \hsize{\hss
    \begin{tabular}{ccccccccccc} \toprule
         Geyser example & \multicolumn{6}{c}{Nb of iterations} & \multicolumn{2}{c}{Mean nb of steps} & Mean & Percent of\\ \cmidrule(lr){2-7} \cmidrule(lr){8-9}
         (dichotomised) & Min & Q1 & Q2 & Mean & Q3 & Max & Fw & Bw & time (s) & convergence\\ \midrule

        qN & 12 & 18 & 21 & 24.23 & 26 & 66 & 24.23 & 0 & 0.33 & 100 \\
        BW & 12 & 82 & 139 & 126 & 161 & 324 & 126 & 126 & 1.21 & 100 \\
        SQUAREM & 15 & 71 & 113 & 108 & 140 & 290 & 78.59 & 71.89 & 0.71 & 100 \\
        QNEM & 6 & 14 & 16 & 16.27 & 19 & 64 & 17.55* & 1.62 & 0.26 & 99.8 \\ \bottomrule
    \end{tabular}\hss}
    \centering
    \caption{This table recapitulates for each algorithm the number of iterations, the mean number of steps, the mean running time and the percentage of convergence in the dichotomised geyser example. The mean running times are in seconds. For all algorithms, each backward step takes $\sim$ 8ms, but the forward step for the Baum-Welch (BW) and the SQUAREM algorithms is $\sim$ 1.5ms only while it is $\sim$ 11.5ms for the quasi-Newton (qN) and the QNEM algorithms, due to the computation of the derivatives. * = The number of forward steps is higher than the number of iterations because the backtracking is not counted in the number of iterations.}
    \label{tab:res gd}
\end{table}

The quasi-Newton and the QNEM algorithms had a lower number of iterations compared to the SQUAREM and the Baum-Welch algorithms. The latters showed a mean number of iterations almost 6 times higher. The lowest number of iterations was always reached by the QNEM. 

Moreover, the quickest algorithm was the QNEM as well, with a mean of 0.26 seconds by run, compared to the quasi-Newton, the Baum-Welch and the SQUAREM algorithms, which had respectively 0.33 seconds, 1.21 seconds and 0.71 seconds. The mean time spent in the forward and backward step is the lower for the QNEM as well, with 215ms compared to 280ms for the quasi-Newton algorithm, 1197ms for the Baum-Welch algorithm and 693ms for the SQUAREM.

We note that the QNEM did not converge for two runs.

\begin{table}[!htbp]
    \hbox to \hsize{\hss
    \begin{tabular}{cccccccccc} \toprule
         Likelihood & \multicolumn{5}{c}{Parameters} & \multicolumn{4}{c}{Percent of runs}\\ \cmidrule(lr){2-6} \cmidrule(lr){7-10}
          & a & b & c & d & e & qN & BW & SQUAREM & QNEM\\ \midrule
         144.5 & 0.79 & 0.57 & 0* & 1* & 0.95 & 52.6 & 50.6 & 50.8 & 42.6 \\ \arrayrulecolor{gray} \midrule
         149.5 & 0.58 & 0.50 & 1* & 0* & 0.57 & 41.2 & 42.4 & 42.4 & 36.2 \\ \midrule
         \multirow{2}{*}{172.8} & 0* & ND & \multirow{2}{*}{0.55 / 0.74} & 0.74 / 0.55 & \multirow{2}{*}{ND} & \multirow{2}{*}{6.2} & \multirow{2}{*}{7.0} & 0.9 & 4.4 \\
          & ND & 0* & & ND & & & & 5.9 & 4.5 \\ \midrule
         Other & & & & & & 0 & 0 & 0 & 11.1 \\ \arrayrulecolor{black} \bottomrule
    \end{tabular}\hss}
    \centering
    \caption{This table shows the different convergence points reached with the dichotomised geyser example and their percent. ND = Non determined. * = These values were rounded for the quasi-Newton results because of the bound-constraints.}
    \label{tab:conv gd}
\end{table}

We then present the table of the convergence points in Table \ref{tab:conv gd}. We see that the optimum is obtained most of the time by all algorithms but it represents only $\sim$ 50$\%$ of the runs for all algorithms except the QNEM, for which it's closer to 40$\%$. There are two other local optimum reached by the algorithms. Moreover, we observe that the QNEM converges to other local optima not reached by the other algorithms.

        \subsubsection{Continuous values}\label{results gc}

For the continuous case, the summary of the results is presented in Table \ref{tab:res gc}.

\begin{table}[!htbp]
    \hbox to \hsize{\hss
    \begin{tabular}{ccccccccccc} \toprule
         Geyser example & \multicolumn{6}{c}{Nb of iterations} & \multicolumn{2}{c}{Mean nb of steps} & Mean & Percent of\\ \cmidrule(lr){2-7} \cmidrule(lr){8-9}
         (continuous) & Min & Q1 & Q2 & Mean & Q3 & Max & Fw & Bw & time (s) & convergence\\ \midrule

        qN & 21 & 60 & 78 & 79.41 & 98 & 159 & 79.41 & 0& 1.47 & 93.2\\
        BW & 8 & 21 & 23 & 25.32 & 26 & 129 & 25.32 & 25.32 & 0.24 & 100\\
        SQUAREM & 12 & 30 & 32 & 34.04 & 35 & 118 & 23.23 & 22.59 & 0.22 & 100\\
        QNEM & 4 & 19 & 23 & 24.01 & 27 & 87 & 31.61* & 1.89 & 0.66 & 99.5\\ \bottomrule
    \end{tabular}\hss}
    \centering
    \caption{This table recapitulates for each algorithm the number of iterations, the mean number of steps, the mean running time and the percentage of convergence in the continuous geyser example. The mean running times are in seconds. For all algorithms, each backward step takes $\sim$ 8ms, but the forward step for the Baum-Welch (BW) and the SQUAREM algorithms is $\sim$ 1.5ms only while it is $\sim$ 17.5ms for the quasi-Newton (qN) and the QNEM algorithms, due to the computation of the derivatives. * = The number of forward steps is higher than the number of iterations because the backtracking is not counted in the number of iterations.}
    \label{tab:res gc}
\end{table}

The quasi-Newton algorithm showed the biggest amount of iterations, with a maximum of 159 iterations. The QNEM keeps showing the less number of iterations globally. 

In term of running time, the lowest are reached by the EM-based algorithms: the Baum-Welch algorithm and the SQUAREM, with a mean running time near 0.23s. Therefore, the quasi-Newton was slow on this continuous case, with a mean running time of 1.47s. The hybrid QNEM showed an intermediate time of 0.66s. This difference of time can be explained by the time spent in the forward step, which is 1.5 millisecond for the forward used in the SQUAREM and the Baum-Welch (conditional probabilities forward) against 17.5 milliseconds for the forward used in the other two algorithms (conditional probabilities forward with derivatives) for a comparable mean number of forward steps.

The quasi-Newton algorithm did not converge for 68 sets of initial parameters and the QNEM for 5 sets of parameters.

\begin{table}[!htbp]
    \hbox to \hsize{\hss
    \begin{tabular}{ccccccccccccc} \toprule
         Likelihood & \multicolumn{8}{c}{Parameters} & \multicolumn{4}{c}{Percent of runs}\\ \cmidrule(lr){2-9} \cmidrule(lr){10-13}
          & a & b & $\mu_s$ & $\mu_l$ & $\mu_{sl}$ & $\sigma_s$ & $\sigma_l$ & $\sigma_{sl}$ & qN & BW & SQUAREM & QNEM\\ \midrule
         265.7 & 0.61 & 0.65 & 2.0 & 4.58 & 4.09 & 0.22 & 0.24 & 0.64 & 19.1 & 40.7 & 40.5 & 35.7 \\ \arrayrulecolor{gray} \midrule
         275.1 & 0.54 & 0.46 & 4.44 & 1.98 & 3.23 & 0.3 & 0.19 & 1.03 & 29.2 & 38.6 & 38.6 & 30.6 \\ \midrule
         303.2 & 0.93 & 0.34 & 2.37 & 1.91 & 4.34 & 0.77 & 0.20 & 0.37 & 21.9 & 9.0 & 9.6 & 9.9\\ \midrule
         316.2 & 1* & 0.28 & 4.45 & ND & 2.79 & 0.30 & ND & 1 & 0.1 & 0 & 0 & 6.4\\ \midrule
         \multirow{2}{*}{335} & 0.44& \multirow{2}{*}{0*} & \multirow{2}{*}{3.72} & 4.24 & 2.0 & \multirow{2}{*}{1.05} & 0.45 & 0.22 & \multirow{2}{*}{7.5} & \multirow{2}{*}{5.5} & \multirow{2}{*}{4.9} & \multirow{2}{*}{5.2} \\ 
          & 0.56 & & & 2.0 & 4.24 & & 0.22 & 0.45 & & & & \\ \midrule
         Other & & & & & & & & & 22.2 & 6.2 & 6.4 & 12.2 \\ \arrayrulecolor{black} \bottomrule
    \end{tabular}\hss}
    \centering
    \caption{This table shows the different convergence points reached with the continuous geyser example and their percent. ND = Non determined. * = These values were rounded for the quasi-Newton results because of the bound-constraints.}
    \label{tab:conv gc}
\end{table}

In Table \ref{tab:conv gc}, the best likelihood is reached most of the time by the Baum-Welch, the SQUAREM and the QNEM algorithms but not by the quasi-Newton algorithm. The latter converges more often to local optima. Most of the time, the quasi-Newton converge to the point with the third best likelihood. This behavior can explain its slow convergence in this example.

    \subsection{HBD segments example}\label{results hbd}

For this last example, the results are presented in Table \ref{tab:res hbd}. 

\begin{table}[!htbp]
    \hbox to \hsize{\hss
    \begin{tabular}{ccccccccccc} \toprule
         HBD segments & \multicolumn{6}{c}{Nb of iterations} & \multicolumn{2}{c}{Mean nb of steps} & Mean & Percent of\\ \cmidrule(lr){2-7} \cmidrule(lr){8-9}
         example & Min & Q1 & Q2 & Mean & Q3 & Max & Fw & Bw & time (s) & convergence\\ \midrule

        qN & 7 & 11 & 14 & 15.36 & 18 & 38 & 15.36 & 0 & 0.32 & 100\\
        BW & 5 & 63 & 73 & 68.53 & 78 & 86 & 68.53 & 68.53 & 2.40 & 100\\
        SQUAREM & 5 & 50 & 57 & 56.19 & 63 & 88 & 38.66 & 38.66 & 1.37 & 100\\
        QNEM & 5 & 9 & 11 & 11.99 & 14 & 27 & 12.83* & 1.40 & 0.32 & 100\\ \bottomrule
    \end{tabular}\hss}
    \centering
    \caption{This table recapitulates for each algorithm the number of iterations, the mean number of steps, the mean running time and the percentage of convergence in the HBD segments example. The mean running times are in seconds. For all algorithms, each backward step takes $\sim$ 30.3ms, but the forward step for the Baum-Welch (BW) and the SQUAREM algorithms is $\sim$ 5.3ms only while it is $\sim$ 20.2ms for the quasi-Newton (qN) and the QNEM algorithms, due to the computation of the derivatives. * = The number of forward steps is higher than the number of iterations because the backtracking is not counted in the number of iterations.}
    \label{tab:res hbd}
\end{table}

The smallest number of iterations is always hit by the QNEM, with a maximum of 27 iterations. It is closely followed by the quasi-Newton algorithm. On the opposite, the Baum-Welch algorithm shows less good performances in term of number of iterations with a mean of 68.53 iterations. The SQUAREM algorithm, is a bit better than the Baum-Welch algorithm with a mean of 56.19 iterations. 

In term of running time, the QNEM and the quasi-Newton shows a similar mean running time of 0.32s, whereas the SQUAREM and the Baum-Welch algorithms are way longer with respectively 1.37s and 2.40s of running time. 
These running times are linked to the time spent in forward and backward steps. Even if the time spent in the forward step by the two EM-based algorithms (5.3ms) is lower than the time spent is the QNEM and quasi-Newton's forward step (around 20ms), the two latter need less forward steps (respectively a mean of 12.83 and 15.36 iterations) than the Baum-Welch and the SQUAREM algorithms (respectively a mean of 68.53 and 38.66 iterations).

For all algorithms, only one likelihood optimum is reach with all 1000 draws of the starting value theta $\theta$. In fact, this likelihood admits only one optimum.

\section{Discussion}

The Baum-Welch (or EM) algorithm is a natural choice when it comes to HMM parameter estimation, offering several advantages: its conceptual simplicity as well as the ease of implementation, and the fact it naturally accommodates the possible constraints on the parameters. However, while the first iterations may move quickly in the space of parameters, once it is close to the solution its convergence is usually slow. 

Another possibility is to use a direct maximisation of the likelihood, particularly via a quasi-Newton method. It is expected to be more efficient than the Baum-Welch algorithm. However, the computation of the likelihood, and of its derivatives, on the logarithmic scale, can be tedious. In the literature, it is usually done using a joint probability version of the forward algorithm (equations \ref{eq:fw a} to \ref{eq:lik}), in which all quantities need to be computed in the logarithmic scale. We showed that it is possible to use the conditional probability version of the forward algorithm (equations \ref{eq:init alpha} to \ref{eq:gamma}) to compute the log-likelihood in a simpler manner. 

A practical problem with the use of the quasi-Newton algorithm is the presence of bounds on the parameters. Typically, in our examples, many parameters take their values in $[0,1]$, and others in $[0, +\infty)$, with possibly an infinite value of the log-likelihood on the border of the space of parameters. The off-the-shelf quasi-Newton method for this case is L-BFGS-B, which tends to favour in its first iterations points on the border of the space, leading to a failure of the algorithm. A common solution, but not totally satisfying, is to use constraints such as $[0.01, 0.99]$ instead.

An appealing alternative is to use an EM acceleration, such as the SQUAREM algorithm. It is simple to implement from an existing Baum-Welch implementation. Bounds on the parameters are easily dealt with, including when the likelihood is infinite on the border.

We propose a new algorithm, QNEM, which alternates between Baum-Welch algorithm iterations and BFGS iterations. Our intuition was that the Baum-Welch algorithm can help to explore efficiently the space of parameters, and the BFGS can allow to achieve a fast convergence once the algorithm is near to a local optimum. The algorithm switches between the two kind of iterations using a convexity criterion. The box-constraints are simply managed with a projection of the gradient on the box when the parameter $\theta$ gets on the border.

We compared the four algorithms on four different examples. The umbrella example (cf. sections \ref{methods um} and \ref{results um}) is a toy example with two hidden states and two possible observed states. As expected, the quasi-Newton was much faster than the Baum-Welch algorithm. The SQUAREM succeeds in accelerating the Baum-Welch algorithm, but is still much slower than the quasi-Newton. Our algorithm, QNEM, is equivalent to the quasi-Newton in term of computational time.

Our second example was based on the Old Faithful data, with three hidden states and two possible observed states (cf. sections \ref{methods gey} and \ref{results gd}). In that case, the QNEM is the faster of all algorithms, followed by the quasi-Newton, then by the SQUAREM, and the Baum-Welch is the slowest.

In our third example, the Old Faithful data were used again, with the same model for the hidden Markov chain, but now the observed states are continuous (cf. sections \ref{methods gey} and \ref{results gc}). Surprisingly, for this model, the quasi-Newton is by far the slowest method. The Baum-Welch algorithm is faster, and the SQUAREM is even slightly faster. The QNEM ranks third, with intermediate performances.

Our last example, which was our motivating example, addresses a concrete question in genetics: making inference about HBD segments on the genome. It has two hidden states, and the observations are discrete (cf. sections \ref{methods hbd} and \ref{results hbd}). In that case, just as in the first example, the QNEM and the quasi-Newton were equivalent, and performed better than the SQUAREM and the Baum-Welch in term of computational time. Our results are in line with those of \citep{bertrand2019rzooroh}, who also compared the L-BFGS-B algorithm to the EM algorithm in term of number of iterations and running time on a more complex HMM model aimed to identify HBD segments. In their work, the quasi-Newton method is several time faster than the EM algorithm.

In all examples except the HBD segments example, there are multiple local likelihood optimum, leading to multiple final estimations of the parameters. In a concrete situation, to deal with this issue, one would need to run the algorithms multiple times from different initialisation of the parameters, to retain the highest likelihood. In our examples, the quasi-Newton and the QNEM seems more prone to converge to a local optimum which is not a global optimum. It is worth to note that in our numeric experiments, the starting points are drawn in very large intervals, which may exacerbate this behaviour.

There was no algorithm uniformly more efficient. In particular, it was surprising to us to find that on our third example, the Baum-Welch and the SQUAREM are the more efficient ones. The presence of multiple optima in this example seems to play a role in this behaviour. The QNEM was the most efficient on our second example, and equivalent to the quasi-Newton on the first and fourth example. This new optimization algorithm for HMM parameter estimation seems to be an interesting alternative to the other methods, offering a good compromise in most situations.

\section*{Code availability}
Related to this paper, we created a R package available at \url{https://github.com/SidonieFoulon/steveHMM}.

\section*{Acknowledgment}

This work is part of the Inserm cross-cutting program GOLD (GenOmics variability in health \& Disease)

\section*{Funding sources}

Sidonie Foulon is funded for this work by a doctoral grant from the French Ministry of Research.

\bibliographystyle{plainnat}
\bibliography{hmmbibli}

\end{document}